\documentclass[]{article}
\usepackage{amsmath,amsthm,amsfonts}

\newtheorem{theorem}{Theorem}[section]
\newtheorem{lemma}[theorem]{Lemma}

\newcommand{\definition}[1]{{\it #1}}

\newcommand{\ses}[3]{
0 \rightarrow #1 \rightarrow #2 \rightarrow #3 \rightarrow 0}

\newcommand{\pushses}[6]{
\begin{matrix}
0& \rightarrow & #1 &\rightarrow & #2 & \rightarrow & #3
&\rightarrow&0\\
 &         &\downarrow &    &\downarrow &      &\downarrow & &  \\
0& \rightarrow & #4 &\rightarrow & #5 & \rightarrow & #6
&\rightarrow&0
\end{matrix}
}
\newcommand{\dual}{^{\vee}}

\DeclareMathOperator{\rk}{rk}
\DeclareMathOperator{\ddim}{\textbf{dim}}
\DeclareMathOperator{\ext}{ext}
\DeclareMathOperator{\Hom}{Hom}
\DeclareMathOperator{\Ext}{Ext}
\DeclareMathOperator{\End}{End}
\DeclareMathOperator{\EEnd}{\mathfrak{End}}

\DeclareMathOperator{\hcf}{hcf}

\begin{document}

\title{Birational classification of moduli spaces of representations
of quivers.}
\author{Aidan Schof\i eld}
\date{November 2, 1999} 
\maketitle

\begin{abstract}
Let $\alpha$ be a Schur root; let $h=\hcf_{v}(\alpha(v))$ and let
$p=1-\langle\alpha/h,\alpha/h\rangle $. Then a moduli space of
representations of dimension vector $\alpha$ is birational to $p$ $h$
by $h$ matrices up to simultaneous conjugacy. Therefore, if
$h=1,2,3\text{ or } 4$, then such a moduli space is a rational variety
and if $h$ divides $420$ it is a stably rational variety. 
\end{abstract}

\section{Introduction}
\label{Intro}

The definitions of terms used without explanation in this section may
be found in the next section.

Suppose that $\alpha$ is a dimension vector for the quiver $Q$. Then
we should like to understand the moduli spaces of the representations
of dimension vector $\alpha$. The representations of dimension vector
$\alpha$ are parametrised by a vector space $R(Q,\alpha)$ on which the
algebraic group $PGl_{\alpha}$ acts so that the orbits correspond to
the isomorphism classes of representations. Thus we are in the general
area of geometric invariant theory. We already know that there is a
stable point for some linearisation of the action of $PGl_{\alpha}$ if
and only if $\alpha$ is a Schur root by theorem 6.1 of \cite{genrep}
and \cite{King} and that all the various moduli spaces obtained are
birational; thus a natural question is to describe the moduli spaces
birationally. This has been attempted in special cases recently
\cite{Dolgachev}, \cite{Megyesi} and \cite{Zaitsev} where the question
attempted has been to show that these spaces are in fact rational
varieties. Another result, in \cite{leBSch}, is that these spaces are
stably birational to a suitable number of suitably sized matrices up
to simultaneous conjugacy. What we shall see in this paper is that
they are actually birational to a suitable number of suitably sized
matrices up to simultaneous conjugacy. This allows rationality to be
proved in all known cases and many more besides; it is also, I would
contend, the correct answer since the rationality of these varieties
of matrices up to simultaneous conjugacy, even their stable
rationality, is problematic.

The proof takes place in several steps. Firstly in section \ref{GKQ}
we deal with the special case of generalised Kronecker quivers; the
$n$th Kronecker quiver has two vertices $v$ and $w$ and $n$ arrows
from $v$ to $w$. Given a dimension vector $\alpha=(a,b)$ where
$\hcf(a,b)=1$, we construct two representations $S$ and $T$ such that
for a general representation $R$ of dimension vector $h\alpha$,
\begin{align*}
\hom(S,R)=&h=\hom(T,R)\\
\intertext{and} 
\ext(S,R)=&0=\ext(T,R).
\end{align*}
Moreover $\hom(S,T)=1+p$ where $p=1-\langle\alpha,\alpha\rangle$. Thus
$\Hom(S\oplus T,\ )$ is a functor which takes most representations of
dimension vector $h\alpha$ to representations of dimension vector
$(h,h)$ for the $(1+p)$th Kronecker quiver. We shall see that this
induces a birational equivalence between their moduli spaces and since
the moduli space of representations of dimension vector $(h\ h)$ for
the $(1+p)$th Kronecker quiver is birational to the moduli space of
$p$ $h$ by $h$ matrices up to simultaneous conjugacy, this completes
the case of Schur roots over a generalised Kronecker quiver.

Next we need a reduction step for an arbitrary Schur root for a
general quiver. This reduction step is defined only for quivers
without loops; however this problem is easily dealt with. In section
\ref{internal}, we shall see that if $\alpha$ is a Schur root for a
quiver without loops then there exist smaller indivisible Schur roots
$\beta$ and $\gamma$ such that a general representation of dimension
vector $\alpha$ has a unique subrepresentation of dimension vector
$c\gamma$ and $\alpha=b\beta+c\gamma$; also
$\hcf_{v}(\alpha(v))=\hcf(b,c)$. Under the inductive assumption that
we have already proved our theorem for all smaller dimension vectors,
this allows a reduction to the case of a dimension vector $(b,c)$ for
a quiver with two vertices $v$ and $w$ whose arrows are either loops
at $v$ and $w$ or arrows from $v$ to $w$.

Finally we handle this case directly in section \ref{finally} using the
results for generalised Kronecker quivers, reduction on the number of
arrows, reflection functors and duality.

The next section summarises the terminology and theory of quivers that
we shall use whilst section \ref{s2} develops the geometric facts that
we shall need. Although the arguments here are rather trivial, they do
appear to be a new tool in the birational study of moduli spaces;
these tools are used in \cite{AKandme} to study moduli spaces of
vector bundles over smooth projective curves and are also an important
ingredient in the proof of a similar result to the present one for the
moduli space of vector bundles over $\mathbb{P}^{2}$.

It seems likely that there should be much better proofs of the main
result of this paper. I believe the following question to have a
positive answer but can prove this only in special cases. Let $Q$ be a
quiver without loops; let $\alpha$ be a Schur root. We say that a
representation $R$ is left perpendicular to the dimension vector
$\alpha$ if for a general representation $S$ of dimension vector
$\alpha$, 
\begin{equation*}
\Hom(R,S)=0=\Ext(R,S).
\end{equation*}
Is there a set of representations
left perpendicular to the dimension vector $\alpha$, $\{R_{i}\}$ such
that $\{R_{i}\}^{\perp}$, the full subcategory of representations $S$
such that 
\begin{equation*}
\Hom(R_{i},S)=0=\Ext(R_{i},S)
\end{equation*}
for all $i$, is equivalent to the
category of finite dimensional representations of a free algebra? This
is a much stronger result if true than the birationality result proved
here. As a less ambitious approach to proving birationality a
different way, it is plausible that the method used to deal with
generalised Kronecker quivers might extend to all quiver without loops
though the role played by preprojective representations would have to
be played by representations $T$ such that $\Ext(T,T)=0$.

\section{Terminology}
\label{s0}

This section contains nothing new; it is a summary of the notation we
shall use and results that are already known.

A \definition{quiver} $Q$ is a directed graph; thus we have a set of
vertices $V$ and a set of arrows $A$ and two functions
$i,t\colon A\rightarrow V$ where we shall write $ia = i(a)$ for the initial
vertex of $a$ and $ta = t(a)$ for the terminal vertex of $a$. A
\definition{representation} $R$ of the quiver $Q$ is a collection of
finite dimensional vector spaces $\{R(v)\colon v\in V\}$ indexed by the
vertex set $V$ and a collection of linear maps $\{R(a):a\in A\}$
indexed by the arrow set where $R(a)\colon R(ia)\rightarrow R(ta)$ is a
linear map from $R(ia)$ to $R(ta)$.  A \definition{homomorphism of representations}
$\phi\colon R\rightarrow S$ is given by a
collection of linear maps $\phi(v)\colon R(v)\rightarrow S(v)$ such that for
all arrows $a$, $R(a)\phi(ta) = \phi(ia)S(a)$. With these definitions
it is a simple matter to see that the category of representations of
the quiver $Q$ is an abelian category.

A \definition{path} $p$ in the quiver of length $n$ is a finite list
of arrows in the quiver $a_{1},\dots,a_{n}$ such that for $j = 1
\text{ to } n-1,\ t(a_{j}) = i(a_{j+1})$. We usually write $p =
a_{1}\dots a_{n}$ and we extend the definition of the functions $i,t$
by $ip = ia_{1}$ and $tp = ta_{n}$. For each vertex $v$ there is a
path of length $0$ at the vertex $v$, $e_{v}$, where we define $ie_{v}
= v = te_{v}$. The \definition{path algebra} $\Pi(Q)$ of the quiver
$Q$ has a basis consisting of all the paths in the quiver. We define a
product in $\Pi(Q)$ by $pq = 0$ if $tp \neq iq$ and $pq = pq$ if $tp =
iq$. One may check that the category of representations of the quiver
$Q$ is the category of finite dimensional modules for the path algebra
$\Pi(Q)$ of the quiver. The path algebra of the quiver is finite
dimensional if and only if there are no loops in the quiver;
that is, there do not exist an integer $n > 0$ and arrows
$a_{1},\dots,a_{n}$ such that for $j = 1$ to $n$, $ta_{j} =
ia_{j+1\mod n}$.

Given a representation $R$ of the quiver $Q$, its
\definition{dimension vector} $\ddim_{Q}(R)$ is a function from the
vertex set to $\mathbb{N}$ defined by $\ddim_{Q}R(v) = \dim R(v)$, the
dimension of $R(v)$ over the field $k$. By \definition{a dimension
vector of the quiver} $Q$ we mean a function from the vertex set to
$\mathbb{N}$. We shall usually use small greek letters from the
beginning of the alphabet for these. We define the
\definition{greatest divisor} $g(\alpha)$ of $\alpha$ to be
$\hcf_{v\in V}(\alpha(v))$. If $g(\alpha)=1$, we say that $\alpha$ is
\definition{indivisible}, and otherwise we say that it is
\definition{divisible}. The \definition{support} of a dimension vector
is the maximal subquiver of $Q$ whose vertex set is $\{v\in
V:\alpha(v)\neq 0\}$. The \definition{vector space of dimension
vectors} is the vector space $\mathbb{R}^{V}$ of functions from $V$ to
$\mathbb{R}$.

We use the notation ${^{a}k^{b}}$ to denote the space of $a$ by $b$
matrices over the field $k$.  Given a dimension vector of the quiver
$Q$, there is a vector space $R(Q,\alpha) = \oplus_{a\in A}
{^{\alpha(ia)}k^{\alpha(ta)}}$ that parametrises representations of
dimension vector $\alpha$. Thus we define for each point $p\in
R(Q,\alpha)$, $R_{p}(v) = k^{\alpha(v)}$ and we define $R_{p}(a)$ by
$p = (\dots,R_{p}(a),\dots)$.  The reductive algebraic group
$Gl_{\alpha}(k)=\times_{v\in V}Gl_{\alpha(v)}(k)$ acts on
$R(Q,\alpha)$ by $R_{pg}(a) = g_{ia}^{-1}R_{p}(a)g_{ta}$ where $g =
(\dots,g_{v},\dots)$. The diagonal embedding of $k^{*}$, the
multiplicative group of the field $k$, acts trivially and thus the
factor group which we call $PGl_{\alpha}(k)=Gl_{\alpha}(k)/k^{*}$ acts
on $R(Q,\alpha)$. We shall usually write $PGl_{\alpha}$ instead of
$PGl_{\alpha}(k)$.  Its orbits are in $1$-to-$1$ correspondence with
the isomorphism classes of representations of dimension vector
$\alpha$. We shall say that a \definition{general} representation has
property $P$ if there exists an open subvariety $U$ of $R(Q,\alpha)$
such that $p\in U\Rightarrow R_{p}$ has property $P$. 

A representation is said to be a \definition{Schur representation} if
its endomorphism ring is $k$. The dimension vector of a Schur
representation is said to be a \definition{Schur root}. The reason
that these are called Schur roots is that the dimension vector of an
indecomposable representation is always a root (in the sense of
Kac-Moody) for a certain symmetric bilinear form on the vector space
of dimension vectors which we shall soon describe.

Given any two representations $R$ and $S$ of the quiver $Q$,
one knows that $\Ext^{2}(R,S)=0$. We use the abbreviations 
\begin{align*}
\hom(R,S)&=\dim \Hom(R,S)\\
\ext(R,S)&= \dim \Ext(R,S)\\
\hom(R,\beta)&= \min_{\ddim S = \beta}\hom(R,S)\\
\hom(\alpha,\beta)&= \min_{\ddim R = \alpha}\hom(R,\beta).
\end{align*}
The notations $\ext(R,\beta)$ and $\ext(\alpha,\beta)$ are defined
analogously. If $\ddim R = \alpha$ and $\ddim S = \beta$ then
\begin{equation*}
\hom(R,S) - \ext(R,S) = \langle \alpha,\beta \rangle =
\sum_{v\in V}\alpha(v)\beta(v) 
- \sum_{a\in A}\alpha(ia)\beta(ta) 
\end{equation*}
which we use as the definition of the \definition{Euler non-symmetric
bilinear form} $\langle \ ,\ \rangle $ on the vector space of
dimension vectors. Kac shows in \cite{Kac1} how to associate to the
symmetrisation of this form a Lie algebra whose set of roots is the
set of dimension vectors of indecomposable representations of the
quiver. When there is no arrow $a$ such that $ia = ta$ then this Lie
algebra is just the Kac-Moody Lie algebra of the underlying graph of
the quiver $Q$ obtained by forgetting the direction of each arrow.
Note that it also follows that
\begin{equation*}
\hom(R,\beta) - \ext(R,\beta) = \hom(\alpha,\beta) -
\ext(\alpha,\beta) = \langle \alpha,\beta \rangle  
\end{equation*}

Given a Schur root $\alpha$, we shall say that it is a
\definition{real} Schur root if $\langle \alpha,\alpha\rangle=1$, that
it is \definition{isotropic} if $\langle \alpha,\alpha\rangle=0$ and
that it is \definition{non-isotropic} in the remaining case where
$\langle \alpha,\alpha\rangle<0$. When $\alpha$ is a real Schur root
then there is a unique indecomposable representation of dimension
vector $\alpha$, $G(\alpha)$. 

We define the \definition{Kac inner product} on the space of dimension
vectors to be the symmetrisation of the Euler non-symmetric bilinear
form, that is, we define $(\alpha,\beta)=\langle \alpha,\beta\rangle +
\langle \beta,\alpha\rangle$. For each vertex $v$, we define a linear
map $r_{v}$ on the vector space of dimension vectors to be reflection
in the hyperplane perpendicular to $e_{v}$ where $e_{v}$ is the
dimension vector $e_{v}(w)=\delta_{vw}$. We say that $\alpha$
\definition{lies in the fundamental region for the action of the Weyl
group} if $r_{v}(\alpha)\geq \alpha$ for all vertices $v$. In the case
where the underlying graph of the quiver is an extended Dynkin diagram
this fundamental region is $1$-dimensional and is the null space of
the Kac form. If $\alpha$ is a dimension vector such that the
underlying graph of its support is an extended Dynkin diagram for
which it lies in the null space we say that it is a \definition{null
dimension vector}. A null dimension vector is a Schur root if and only
if it is indivisible. We note that Kac shows in \cite{Kac1} that all
dimension vectors with connected support in the fundamental region for
the action of the Weyl group that are not divisible null roots are in
fact Schur roots. Thus if $\alpha$ is a dimension vector such that
$r_{v}(\alpha)\geq \alpha$ for all vertices $v$, then $\alpha$ is a
Schur root if and only if its support is connected and it is not a
divisible null dimension vector.

Given a quiver $Q$, we say that a vertex $v$ is a \definition{source}
if it is not the terminal vertex of any arrow. Similarly we say that
it is a \definition{sink} if it is not the initial vertex of any
arrow.  If $v$ is a source or a sink, there is a quiver $Q_{v}^{+}$
when it is a source or a quiver $Q_{v}^{-}$ when it is a sink, with
the same vertex and arrow set as $Q$; however, the incidence functions
on $Q_{v}^{+}$ or $Q_{v}^{-}$, $i_{v}$ and $t_{v}$ from $A$ to $V$ are
defined by $i_{v}a=ia$ and $t_{v}a=ta$ if $a$ is not incident to $v$
and $i_{v}a=ta,\ t_{v}a=ia$ if $a$ is incident to $v$. If $v$ is a
source, and $R$ is a representation of the quiver $Q$, we define a
vector space $e_{v}^{+}(R)(v)$ as the cokernel of the homomorphism
$\oplus_{a,ia=v}R(a)$ from $R(v)$ to $\oplus_{a,ia=v}R(ta)$ and if $v$
is a sink, then we define a vector space $e_{v}^{-}(R)(v)$ to be the
kernel of the homomorphism $\oplus_{a,ta=v}R(a)$ from
$\oplus_{a,ta=v}R(ia)$ to $R(v)$. When $v$ is a source there is a
functor $e_{v}^{+}$ from the category of representations of the quiver
$Q$ to the category of representations of the quiver $Q_{v}$ defined
by $e_{v}^{+}(R)(w)=R(w)$ for $w\neq v$ and $e_{v}^{+}(R)(a)=R(a)$ for
$a$ not incident to $v$ whilst for $a$ incident to $v$ we define
$e_{v}^{+}(R)(a)$ to be the $a$-th component of the homomorphism from
$\oplus_{a,ia=v}R(ta)$ to $e_{v}^{+}(R)(v)$. When $v$ is a sink the
functor $e_{v}^{-}$ is defined in a similar way. $e_{v}^{+}$ and
$e_{v}^{-}$ are called \definition{reflection functors} at the vertex
$v$ and are mutually inverse between the full subcategories of
representations that have no summand isomorphic to the simple
representation at the vertex $v$. If the representation $R$ has no
simple summand at the vertex $v$ then
$\ddim_{Q_{v}^{*}}(e_{v}^{*}(R))=r_{v}(\ddim_{Q}(R))$ for $*=+$ or
$-$. If $v$ is a source for the quiver $Q$ and $\alpha$ is a dimension
vector such that there are representations of dimension vector
$\alpha$ having no summands isomorphic to the simple at the vertex
$v$, or equivalently, $\alpha(v)\leq \oplus_{a,ia=v}\alpha(ia)$ then
the reflection functors allow us to show that a moduli space of
representations of the quiver $Q$ of dimension vector $\alpha$ is
birational to a moduli space of representations of the quiver
$Q_{v}^{+}$ of dimension vector $r_{v}(\alpha)$. 

Let us fix a dimension vector $\alpha$. A representation $R$ of
dimension vector $\alpha$ is isomorphic to a direct sum of
indecomposable representations $R\cong \oplus_{i}R_{i}$ and this
determines a sum $\alpha=\sum_{i}\beta_{i}$ where $\beta_{i}$ is the
dimension vector of $R_{i}$. Kac notes in \cite{Kac1} and \cite{Kac2}
that this sum is constant on constructible pieces of $R(Q,\alpha)$ and
hence there exists a dense open subvariety $U$ of $R(Q,\alpha)$ and a
particular sum $\alpha=\sum_{i}\beta_{i}$ known as the
\definition{canonical decomposition} of $\alpha$ such that for every
$p\in U$, $R_{p}\cong \oplus_{i}R_{p,i}$ where each $R_{p,i}$ is an
indecomposable representation. In fact, Kac shows that each
$\beta_{i}$ in the canonical decomposition of $\alpha$ is a Schur root
and a general representation of dimension vector $\alpha$ is
isomorphic to a direct sum $\oplus_{i}R_{i}$ where for each $i$,
$R_{i}$ is a representation of dimension vector $\beta_{i}$ such that
$\End(R_{i})=k$, the ground field. We call each $\beta_{i}$ a
\definition{canonical summand} of $\alpha$.

There is another number of importance to us in the case where $\alpha$
is a Schur root. We define the \definition{parameter number} of
$\alpha$, $p(\alpha)$, to be $1-\langle \alpha,\alpha \rangle$. A
dimension count shows that this is the expected dimension of a moduli
space of representations of dimension vector $\alpha$. Let $h$ be the
greatest divisor of the Schur root $\alpha$ and let $p=p(\alpha/h)$;
we shall see that a moduli space of representations of dimension
vector $\alpha$ is birational to $p$ $h$ by $h$ matrices up to
simultaneous conjugacy.

\section{Birational geometry}
\label{s2}

Let $X$ be an algebraic variety on which the algebraic group $G$
acts. Let
\begin{equation*}
1\rightarrow k^{*}\rightarrow \tilde{G}\rightarrow G\rightarrow 1
\end{equation*}
be a short exact sequence of algebraic groups.  Let $E$ be a vector
bundle over $X$ on which $\tilde{G}$ acts compatibly with the action
of $G$. Then $k^{*}$ acts on the fibres of $E$ and if this action is
via the character $\phi_{w}(\lambda)=\lambda^{w}$ then we shall say
that $E$ is \definition{a $\tilde{G}$ vector bundle of weight} $w$. A
morphism of $\tilde{G}$ vector bundles of weight $w$ is a morphism of
vector bundles that is also $\tilde{G}$ equivariant. 

A \definition{family of representations} $\mathcal{R}$ of the quiver
$Q$ over the algebraic variety $X$ is a collection of vector bundles
$\{\mathcal{R}(v):v\in V\}$ and morphisms of vector bundles
$\{\mathcal{R}(a):a\in A\}$ where $\mathcal{R}(a)\colon
\mathcal{R}(ia)\rightarrow\mathcal{R}(ta)$. For each point $p\in X$,
there is a representation $\mathcal{R}_{p}$ of the quiver. We shall
say that the family is \definition{general} if there exists an open
subvariety $U$ of $R(Q,\alpha)$ such that for each point $q\in U$
there exists a point $p\in X$ such that $R_{q}\cong
\mathcal{R}_{p}$. We shall say that the family is
\definition{$G$-general} if it is a general family, $X$ is affine and
the algebraic group $G$ acts on $X$ freely so that
$\mathcal{R}_{p}\cong \mathcal{R}_{q}$ if and only if $p$ and $q$ lie
in the same orbit. Finally we say that the family is
\definition{$(\tilde{G},G)$-standard} or that it is a
\definition{standard $(\tilde{G},G)$ family} if and only if it is
$G$-general and each $\mathcal{R}(v)$ is a $\tilde{G}$ vector bundle
of weight $1$. The assumption that $X$ is affine is simply to avoid
problems about the existence of certain orbit space.

Let $\alpha$ be some dimension vector and take $G$ to be
$PGl_{\alpha}$ whilst $\tilde{G}$ is $Gl_{\alpha}$. Then the usual
family on $R(Q,\alpha)$ satisfies this last condition. If $\alpha$ is
in addition a Schur root then there is an equivariant affine open
subvariety of $R(Q,\alpha)$ consisting entirely of stable points for
some suitable linearisation of the action of $PGl_{\alpha}$ and the
restriction of the usual family on $R(Q,\alpha)$ is a standard
$(Gl_{\alpha},PGl_{\alpha})$ family. If $\alpha$ is a Schur root,
$h=\hcf_{v}(\alpha(v))$ and $p=1-\langle \alpha/h,\alpha/h\rangle $,
we shall say that $\alpha$ is \definition{reducible to matrix normal
form} if there exists a $(Gl_{h},PGl_{h})$-standard family of
representations of dimension vector $\alpha$, $\mathcal{R}$, over an
algebraic variety $X$ such that $X$ is $PGl_{h}$-birational to
$M_{h}(k)^{p}$ where $PGl_{h}$ acts by conjugation on each factor of
$M_{h}(k)^{p}$. 

We shall need to consider some generalities on $\tilde{G}$ vector
bundles of weight $w$. It is important to show that two such vector
bundles of the same rank are locally $\tilde{G}$ isomorphic.  Let $E$
be a vector bundle over the affine algebraic variety $X$ on which $G$
acts freely.  If $E$ is a $\tilde{G}$ vector bundle of weight $w$ then
$E\dual$ is a $\tilde{G}$ vector bundle of weight $-w$; if in addition
$F$ is a $\tilde{G}$ vector bundle of weight $w'$ then $E\otimes F$ is
a $\tilde{G}$ vector bundle of weight $w+w'$. Thus $E\dual\otimes E$
is a $\tilde{G}$ vector bundle of weight $0$, that is, $G$ acts on
$E\dual\otimes E$ compatibly with the action of $G$ on $X$ and
therefore freely. Therefore, $E\dual\otimes E/G$ is a vector bundle
over $X/G$ and in fact it is a bundle of central simple algebras over
$X/G$ since $G$ acts as automorphisms of the sheaf of algebras
$E\dual\otimes E=\EEnd(E)$. Thus the fibre of $E\dual\otimes E/G$ over
the generic point of $X/G$ is a central simple algebra over the
function field of $X/G$ and it is isomorphic to $M_{n}(D)$ for a
suitable central division algebra $D$. We can read off the dimension
of $D$ from the ranks of $\tilde{G}$ vector bundles of weight $w$ over
$G$ equivariant open subvarieties of $X$ and this in turn allows us to
determine the local isomorphism of such $\tilde{G}$ vector bundles of
the same rank and weight.

\begin{lemma}
\label{lem:tg}
Let $X$ be an affine algebraic variety on which the algebraic group
$G$ acts and let $E$ and $F$ be $\tilde{G}$ vector bundles of weight
$w$. Then the sheaves of algebras over $X/G$, $E\dual\otimes E/G$ and
$F\dual\otimes F/G$ are Morita equivalent. Let $D_{w}$ be a central
division algebra over the function field $\mathfrak{F}(X/G)$ of the
generic point of $X/G$ that is Morita equivalent to the fibre of
$E\dual\otimes E/G$ over the generic point of $X/G$. Then the
dimension of $D_{w}$ over $\mathfrak{F}(X/G)$ is $n^{2}$ where $n$ is
the minimal rank of a vector bundle of weight $w$ over some
equivariant open subvariety of $X$ or equivalently the highest common
factor of the ranks of vector bundles of weight $w$ over equivariant
open subvarieties of $X$. Further, if $\rk(E)\leq \rk(F)$, then there
exists an equivariant open subvariety $Y$ of $X$ such that the
restriction of $E$ to $Y$ is $\tilde{G}$ isomorphic to a $\tilde{G}$
summand of the restriction of $F$ to $Y$. In particular, if
$\rk(E)=\rk(F)$, then $E$ and $F$ are isomorphic over $Y$.
\end{lemma}
\begin{proof}
The vector bundle $E\dual\otimes F$ is a vector bundle of weight $0$
so we may form the vector bundle over $X/G$, $E\dual\otimes F/G$, and
the sheaf of algebras $E\dual\otimes E/G$ acts on its left whilst
$F\dual\otimes F/G$ acts on the right. Moreover, it is clear that 
$(E\dual\otimes E/G)^{o}\otimes (F\dual\otimes F/G)$ is isomorphic to
$\EEnd(E\dual\otimes F/G)$ where $(E\dual\otimes E/G)^{o}$ is the sheaf
of opposite algebras to the sheaf of algebras $(E\dual\otimes E/G)$ and thus 
our two sheaf of algebras are Morita equivalent as stated. 

Thus if $n^{2}$ is the dimension of $D_{w}$ over $\mathfrak{F}(X/G)$,
and the rank of $E$ is $e$, we see that $n^{2}|e^{2}$ and hence $n|e$;
let $t=e/n$.  Thus $n$ divides the highest common factor of the ranks
of $\tilde{G}$ vector bundles of weight $w$ over $G$ equivariant open
subvarieties of $X$. The fibre of $E\dual\otimes E/G$ over the
generic point of $X/G$ is isomorphic to $M_{t}(D_{w})$. Therefore,
there exists an affine open subvariety $Y/G$ of $X/G$ on which
$E\dual\otimes E/G\cong M_{t}(A)$ for a suitable sheaf of algebras
$A$ on $Y/G$. The matrix units in this sheaf of algebras which are
$\tilde{G}$ invariant endomorphisms of $E$ show that $E\cong
E_{1}^{t}$ where $E_{1}$ is a $\tilde{G}$ vector bundle of weight $w$
on $Y$ whose rank is $n$. It follows that $n$ is actually the minimal
rank of a $\tilde{G}$ vector bundle of weight $w$ over a $G$
equivariant open subvariety of $X$ and since it also divides all these
ranks it is both the minimum and the highest common factor. 

For the final part, consider $\EEnd(E\oplus F)/G$. The fibre above the
generic point of $X/G$ is $M_{t+u}(D_{w})$ where $\rk(F)=nu$ and all
idempotents of the same rank are conjugate, thus there exists an
affine $G$ equivariant open subvariety $Y$ on which there is a
$\tilde{G}$ vector bundle $E_{1}$ of weight $w$ and rank $n$ such that
$E$ is $\tilde{G}$ equivariantly isomorphic to $E_{1}^{t}$ and $F$ is
$\tilde{G}$ equivariantly isomorphic to $E_{1}^{u}$ which implies the
last sentence of the lemma.
\end{proof}

We shall need results that allow us to deduce that a dimension vector
for a quiver is reducible to matrix normal form provided that a
certain related dimension vector for a related quiver is reducible to
matrix normal form. We end this section with three such results; the
first allows us to remove an arrow provided that the dimension vector
remains a Schur root for the smaller quiver; the second shows that
being reducible to matrix normal form is well-behaved with respect to
reflection functors; the third shows that it is also well-behaved
under passage to the dual quiver. If $Q$ is a quiver then the
\definition{dual quiver} $Q\dual$ has the same vertex set as $Q$ and
there is a bijection $\dual{}\colon A\rightarrow A\dual$ where $A\dual$
is the arrow set of $Q\dual$ and $ia\dual =ta$ and
$ta\dual =ia$. Vector space duality gives a contravariant functor from
the category of representations of the quiver $Q$ to the category of
representations of the quiver $Q\dual$.

\begin{lemma}
\label{addanarrow}
Suppose that the dimension vector $\alpha$ for the quiver $Q$ is
reducible to matrix normal form. Let $Q'$ be the quiver obtained by
adjoining one arrow from the vertex $w_{1}$ to the vertex
$w_{2}$. Then $\alpha$ is reducible to matrix normal form over the
quiver $Q'$.
\end{lemma}
\begin{proof}
Let $b$ be the new arrow from $w_{1}$ to $w_{2}$ and suppose that
$\alpha(v)=hc$ and $\alpha(w)=hd$. Let $\mathcal{R}$ be a
$(Gl_{h},PGl_{h})$-standard family of representations over $X$ of
dimension vector $\alpha$ for the quiver $Q$. Then the vector bundle
$\mathcal{E}=\mathcal{R}(w_{1})\dual\otimes \mathcal{R}(w_{2})$ over
$X$ carries a family $\mathcal{R}'$ of representations of dimension
vector $\alpha$ over the quiver $Q'$. Here
$\mathcal{R}'(v)=\mathcal{R}(v)$ and $\mathcal{R}'(a)=\mathcal{R}(a)$
for $a\in Q$ and $\mathcal{R}'(b)$ is the universal map from
$\mathcal{R}(w_{1})$ to $\mathcal{R}(w_{2})$.  Moreover the orbits of
$PGl_{h}$ correspond to the isomorphism classes of representations of
$Q'$ since $PGl_{h}$ acts freely on $X$. The vector bundle
$\mathcal{E}$ is a bundle of weight $0$ and hence by lemma
\ref{lem:tg} there exists an open subvariety of $X$ on which
$\mathcal{E}$ is $PGl_{h}$-isomorphic to
$\mathcal{F}=M_{h}(k)^{cd}\times X$ where $PGl_{h}$ acts by
conjugation on $M_{h}(k)$ and if $X$ is an open affine subvariety of
$M_{h}(k)^{n}$ for some integer $n$ then $\mathcal{F}$ is an open
affine subvariety of $M_{h}(k)^{n+cd}$.
\end{proof}

Next we should deal with reflection functors. Let $v$ be a source for
the quiver $Q$ and let $\mathcal{R}$ be a family over the algebraic
variety $X$ of representations of the quiver $Q$. If there is a point
$p$ such that the simple representation $S_{v}$ is not a summand of
$\mathcal{R}_{p}$, then there is a non-empty open subvariety of $X$ where
this holds since having such a summand is a closed condition and if
$X$ is $G$-general we may choose this subvariety to be $G$-general as
well. Thus we shall assume that for all $p\in X$, $\mathcal{R}_{p}$
has no summand isomorphic to $S_{v}$. Then  
$\phi\colon\mathcal{R}(v)\rightarrow \oplus_{a,ia=v}\mathcal{R}(ta)$
is an injective morphism of vector bundles and we define its cokernel
to be $e_{v}^{+}\mathcal{R}(v)$. For all vertices $w\neq v$, we define
$e_{v}^{+}\mathcal{R}(w)=\mathcal{R}(w)$. It is clear how to define
$e_{v}^{+}\mathcal{R}(a)$ for each arrow $a$ of $Q_{v}^{+}$.The next
lemma shows that the family of representations of the quiver
$Q_{v}^{+}$, $e_{v}^{+}\mathcal{R}$ over $X$ inherits most of the
properties of $\mathcal{R}$. 

\begin{lemma}
\label{lem:reflection}
Let $v$ be a source for the quiver $Q$. Let $\mathcal{R}$ be a family
over the algebraic variety $X$ of representations of the quiver such
that for each point $p$ of $X$, the simple representation at the
vertex $v$ is not a summand of $\mathcal{R}_{p}$. The resulting
family of representations of the quiver $Q_{v}^{+}$,
$e_{v}^{+}\mathcal{R}$, is general, $G$-general, or
$(\tilde{G},G)$-standard if and only if $\mathcal{R}$ is.  In
particular, a dimension vector $\alpha$ such that
$\sum_{a,ia=v}\alpha(ta)\geq \alpha(v)$ for the quiver $Q$ is
reducible to matrix normal form if and only if $e_{v}(\alpha)$ for the
quiver $Q_{v}^{+}$ is reducible to matrix normal form.
\end{lemma}
\begin{proof}
Take a point $p\in X$ and let $U$ be an open subvariety of
$R(Q,\alpha)$ containing a point $q$ such that $R_{q}\cong
\mathcal{R}_{p}$ and for all points $u\in U$, there exists a point
$x\in X$ such that $R_{u}\cong \mathcal{R}_{x}$. Let $V$ be an open
subvariety of $R(Q_{v}^{+},e_{v}^{+}(\alpha))$ such that there exists
a point $q'\in V$ where $R_{q'}\cong e_{v}^{+}(R_{q})\cong
e_{v}^{+}(\mathcal{R}_{p})$, and for all $v\in V$, $R_{v}$ has no
summand isomorphic to $S_{v}$, the simple representation at the vertex
$v$ and $e_{v}^{-}(R_{v})$ is isomorphic to $R_{u}$ for some $u\in U$.
This is possible since the last two conditions are open and must hold
in some neighbourhood of any point $q'$ such that $R_{q'}\cong
e_{v}^{+}(\mathcal{R}_{p})$. Then if $R_{u}\cong \mathcal{R}_{x}$,
$R_{v}\cong e_{v}^{+}(\mathcal{R}_{x})\cong
e_{v}^{+}(\mathcal{R})_{x}$ which shows that $e_{v}^{+}(\mathcal{R})$
is general if $\mathcal{R}$ is. That it is $G$-general when
$\mathcal{R}$ is follows at once and that it is
$(\tilde{G},G)$-standard when $\mathcal{R}$ is follows from the fact
that $e_{v}^{+}(\mathcal{R})(v)$ is a vector bundle of weight equal to
that of each $\mathcal{R}(w)$.
\end{proof}

Let $\mathcal{R}$ be a family of representations of the quiver $Q$; we
define a family of representations of the quiver $Q\dual$,
$\mathcal{R}\dual$, by $\mathcal{R}\dual (v)=\mathcal{R}(v)\dual$ and
$\mathcal{R}(a\dual)=\mathcal{R}(a)\dual$. If $\mathcal{R}$ is a
$(\tilde{G},G)$-standard family, then the natural action of
$\tilde{G}$ on $\mathcal{R}\dual$ means that each $\mathcal{R}\dual
(v)$ is a vector bundle of weight $-1$ so we shall assume that there
is an automorphism $\sigma$ of $\tilde{G}$ such that
$\sigma(\lambda)=\lambda^{-1}$ for $\lambda\in k^{*}$ and we shall
regard $\mathcal{R}\dual$ as acted on by $\tilde{G}$ via this
automorphism. In the case where $G$ is $PGl_{\alpha}$ and $\tilde{G}$
is $Gl_{\alpha}$ then the transpose inverse automorphism on each
factor $Gl_{\alpha(v)}$ of $Gl_{\alpha}$ is what we need.

\begin{lemma}
\label{lem:dual}
Let $\mathcal{R}$ be a family of representations of the quiver
$Q$. Assume that $\tilde{G}$ has an automorphism $\sigma$ such that
$\sigma(\lambda)=\lambda^{-1}$ for $\lambda\in k^{*}$. 
Then $\mathcal{R}\dual$ is general, $G$-general or
$(\tilde{G},G)$-standard if and only if $\mathcal{R}$ is. Let $\alpha$
be a dimension vector for the quiver $Q$. Then $\alpha$ as a dimension
vector for the quiver $Q$ is reducible to matrix normal form if and
only if $\alpha$ as a dimension vector for the quiver $Q\dual$ is
reducible to matrix normal form.
\end{lemma}
\begin{proof}
Taking the same precautions as in the proof of the preceding lemma
proves this result in the same way.
\end{proof}

\section{Generalized Kronecker quivers}
\label{GKQ}

In this section, we shall show that a Schur root for a generalised
Kronecker quiver is reducible to matrix normal form. We shall show
this by reduction to the case of the dimension vector $(h\ h)$ for the
generalised Kronecker quiver $Q(n)$. Thus we should first of all show
that this dimension vector is reducible to matrix normal form. This
may be done as follows. Consider the algebraic variety
$X=M_{h}(k)^{n-1}$; this has $n-1$ canonical morphisms of vector
bundles $\phi_{i}$ for $i=1\rightarrow n-1$ from $\mathcal{O}_{X}^{h}$
to itself given by the $n-1$ components of $X$. We define a family of
representations of dimension vector $(h\ h)$ for the quiver $Q(n)$
over $X$ by defining
$\mathcal{R}(v)=\mathcal{O}_{X}^{h}=\mathcal{R}(w)$,
$\mathcal{R}(a_{0})=I_{\mathcal{O}_{X}^{h}}$ and
$\mathcal{R}(a_{i})=\phi_{i}$ for $0<i<n$. This is a general family
since for a general representation $R$ of dimension vector $(h\ h)$
for the quiver $Q(n)$, $R(a_{0})$ is invertible and therefore is
isomorphic to $\mathcal{R}_{p}$ where the $i$th component of $p$ is
$R(a_{0})^{-1}R(a_{i})$. Moreover, $PGl_{h}$ acts on $X$ so that
$\mathcal{R}_{p}\cong \mathcal{R}_{q}$ if and only if $p$ and $q$ lie
in the same orbit for $PGl_{h}$. Certainly $\mathcal{R}(v)$ and
$\mathcal{R}(w)$ are vector bundles of weight $1$ for $Gl_{h}$;
however, it is clear that $PGl_{h}$ does not act freely. However, the
subvariety of points with non-trivial stabiliser is closed and
$PGl_{h}$ equivariant and there exists a $PGl_{h}$ invariant function
that vanishes on this closed subvariety, so the restriction of the
family $\mathcal{R}$ to the affine open subvariety where this
invariant function is non-zero is a $(Gl_{h},PGl_{h})$ standard family
as required.

Next we note that the study of general representations for
generalised Kronecker quivers is rather simpler than for arbitrary
quivers. We recall certain special representations of the quiver
$Q(n)$. There is a unique source vertex and a unique sink vertex and thus
it makes sense to talk of $e^{+}$ and $e^{-}$ instead of $e_{v}^{+}$
and $e_{w}^{-}$. Moreover, $Q(n)^{+}$ and $Q(n)^{-}$ are isomorphic to
$Q(n)$. The simple representation $S_{w}$ has dimension vector $(0\
1)$ and is projective; we shall call it P(0). By induction, we define
the $n$th \definition{preprojective representation} $P(n)$ to be
$e^{+}P(n-1)$. Also $Q(n)\dual$ is isomorphic to $Q(n)$ and we may
therefore define the \definition{preinjective representations} of
$Q(n)$ by $I(n)=P(n)\dual$. Note that $P(1)$ is the other projective
representation. A representation of the type $P(i)^{c}\oplus
P(i+1)^{d}$ will be called a \definition{general preprojective
representation}.

\begin{theorem}
\label{KQ1}
Let $(a\ b)$ be a dimension vector for the $n$th Kronecker
quiver. Then either there exists an integer n such that either
${e^{-}}^{n}(a\ b)$ or ${e^{+}}^{n}(a\ b)$ is a root in the
fundamental region for the action of the Weyl group and consequently
$(a\ b)$ is a Schur root or else a general representation of dimension
vector $(a\ b)$ is of the form $P(l)^{c}\oplus P(l+1)^{d}$ if $a<b$
and of the form $I(l)^{c}\oplus I(l+1)^{d}$ if $b<a$. In particular,
all roots are Schur roots.
\end{theorem}
\begin{proof}
By lemma \ref{lem:dual} and the above discussion we may assume that
$a\leq b$. If $na\leq b$, then a general representation of dimension
vector $(a\ b)$ is isomorphic to $P(1)^{a}\oplus P(0)^{b-na}$. So we
shall assume that $na>b$. If $b\leq na-b$ then $(a\ b)$ lies in the
fundamental region for the action of the Weyl group. If neither of
these occur then $e^{-}(a\ b)=(a\ na-b)$ is a smaller dimension vector
and if $na-b>a$ then this lies in the fundamental region for the
action of the Weyl group. Otherwise, after dualising (see lemma
\ref{lem:dual}) we have reached the dimension vector $(na-b\ a)$ for
the quiver $Q(n)$. In either case the result follows by induction on
the dimension vector.
\end{proof}

Thus in studying moduli spaces of representations of dimension vector
the Schur root $(a\ b)$ either it is a real root in which case, there
is only one indecomposable representation and our result is trivially
true or else a series of reflection functors lead to a root in the
fundamental region and it is enough to deal with these by
lemma \ref{lem:reflection}.  

We shall say that a dimension vector $(a\ b)$ is a
\definition{pre-projective} dimension vector if a general
representation of dimension vector $(a\ b)$ is isomorphic to
$P(l)^{c}\oplus P(l+1)^{d}$ and that it is a
\definition{pre-injective} dimension vector if a general
representation of dimension vector $(a\ b)$ is isomorphic to
$I(l)^{c}\oplus I(l+1)^{d}$; thus every dimension vector is
pre-projective, pre-injective or a Schur root.

\begin{lemma}
\label{lem:prep}
Let $(a\ b)$ be a root in the fundamental region of the Weyl group for
the $n$th Kronecker quiver such that $\hcf(a,b)=1$ and $a\leq b$. Let
$(c\ d)$ be the dimension vector such that $\langle
(c\ d),(a\ b)\rangle=1$ and $a<c\leq a+b$. Let $S$ be a general
representation of dimension vector $(c\ d)$. Then for a general
representation $R$ of dimension vector $(a\ b)$, $\hom(S,R)=1$,
$\ext(S,R)=0$ and the kernel of the non-zero homomorphism from $S$ to
$R$ is a general preprojective representation $K$.
\end{lemma}
\begin{proof}
We begin by showing that $(c-a\ d-b)$ is a preprojective dimension
vector. 

Firstly, $\langle (c\ d),(a\ b)\rangle=db-c(nb-a)$ and so there is a
unique dimension vector $(c\ d)$ such that $\langle (c\ d),(a\
b)\rangle=1$ and $a<c\leq a+b$; in particular, $c-a\leq b$ .  Since
$\langle (c\ d),(a\ b)\rangle=1$, it follows that
\begin{equation*}
(d-b)b-(c-a)(nb-a)=1+nab-a^{2}-b^{2} 
\end{equation*}
and so 
\begin{align*}
\frac{d-b}{c-a}&=n-\frac{a}{b}+\frac{1+nab-a^{2}-b^{2}}{(c-a)b} \\
               &\geq n-\frac{a}{b}+\frac{1+nab-a^{2}-b^{2}}{b^{2}}\\
	       &> n+\frac{(n-1)ab-a^{2}-b^{2}}{b^{2}}\\
               &= n - (\lambda^{2}-(n-1)\lambda+1) 
\end{align*}
where $\lambda=a/b$. Since $(a\ b)$ is in the fundamental region,
$2/n\leq \lambda=a/b\leq 1$ and the maximal value of
$\lambda^{2}-(n-1)\lambda+1$ on this interval is for $\lambda=2/n$.
Therefore, 
\begin{equation*}
\frac{d-b}{c-a}>n-\frac{-n^{2}+2n+4}{n^{2}}. 
\end{equation*}
Thus for $n\geq 4$, $\tfrac{d-b}{c-a}>n$ and for $n=3$,
$\tfrac{d-b}{c-a}>3-\tfrac{1}{9}>\tfrac{8}{3}$. Thus when $n\geq 4$,
the dimension vector $(c-a\ d-b)$ is the dimension vector of a
projective representation and when $n=3$ it is either the dimension
vector of a projective representation or else of a representation of
type $P(1)^{e}\oplus P(2)^{f}$ (note that the dimension vector of
$P(2)$ is $(3\ 8)$). In either case it is a preprojective dimension
vector. 

Next we see that $\ext((c\ d),(a\ b))=0$. To calculate this, we need
to show that if $(e\ f)$ is a dimension vector of a subrepresentation
of a general representation of dimension vector $(c\ d)$, then
$\langle (e\ f),(a\ b)\rangle\geq 0$. If $(c\ d)$ is a preprojective
dimension vector then this is clear since the possibilities for $(e\
f)$ are all themselves preprojective dimension vectors. If $(c\ d)$ is
not a preprojective dimension vector then it must be a Schur root by
theorem \ref{KQ1} and hence there are stable representations and this
implies that when $(e\ f)$ is a dimension vector of a
subrepresentation of a general representation of dimension vector $(c\
d)$, then $e/f<c/d$ and since $\langle (c\ d),(a\ b)\rangle =0$, it
follows that $\langle (e\ f),(a\ b)\rangle >0$.

Therefore by section 5 of \cite{genrep}, if $R$ is a general
representation of dimension vector $(a\ b)$ and $S$ is a general
representation of dimension vector $(c\ d)$ then $\hom(S,R)=1$, the
unique homomorphism is surjective and the kernel is general of
dimension vector $(c-a\ d-b)$ and hence must be general preprojective
as required.
\end{proof}

Let 
\begin{equation*}
\ses{K}{S}{R}
\end{equation*}
 be a short exact sequence as constructed in the
previous lemma. Applying $\Hom(K,\ )$ and noting that $\Ext(K,K)=0$
shows that the map from $\Hom(K,S)$ to $\Hom(K,R)$ is surjective. Now
we apply $\Hom(\ ,R)$ and note that $\Ext(S,R)=0$ to deduce that the
natural map from $\Hom(K,R)$ to $\Ext(R,R)$ is also surjective. Let
$W$ be a vector subspace of $\Hom(K,S)$ mapped isomorphically to
$\Ext(R,R)$ by the composition of these two surjective maps. There is
a natural map $\phi$ from $K$ to $W\dual\otimes S$.  We define a new
representation $T$ by the folowing pushout diagram:
\begin{equation*}
\pushses{K}{S}{R}{W\dual\otimes S}{T}{R}. 
\end{equation*}

\begin{lemma}
\label{lem:tilt}
Let $K$, $R$, $S$ and $T$ be representations as constructed
above. Then 
\begin{enumerate}
\item $\hom(T,R)=1$, $\ext(T,R)=0$ and the homomorphism from $T$ to
$R$ is surjective.
\item $\hom(S,T)=1+p$ where $p=1-\langle (a\ b),(a\ b)\rangle$.
\item Let $R'$ be a general representation of dimension vector $h(a\
b)$; then 
\begin{equation*}
\hom(T,R')=h=\hom(S,R').
\end{equation*}
Let $W'$ be the kernel of the linear map 
from $\Hom(S,T)\otimes \Hom(T,R')$ to $\Hom(S,R')$. Then the complex
below is a short exact sequence.
\begin{equation*}
\ses{W'\otimes S}{\Hom(T,R')\otimes T}{R'}. 
\end{equation*}
\end{enumerate}
\end{lemma}
\begin{proof}
Apply $\Hom(\ ,R)$ to the short exact sequence 
\begin{equation*}
\ses{W\dual\otimes S}{T}{R}.
\end{equation*}
Then by construction the homomorphism from $W\cong
\Hom(W\dual\otimes S,R)$ to $\Ext(R,R)$ is an isomorphism. Therefore
the first part of this lemma follows.

Apply $\Hom(S,\ )$ to the short exact sequence 
\begin{equation*}
\ses{W\dual\otimes S}{T}{R}.
\end{equation*}
By construction, the non-zero homomorphism from $S$ to $R$
lifts through $T$ and hence the second part of this lemma follows.

Since 
\begin{equation*}
\ext(T,R^{h})=0=\ext(S,R^{h})
\end{equation*}
it follows that for a general
representation of dimension vector $h(a\ b)$, $R'$,
\begin{equation*}
\ext(T,R')=0=\ext(S,R')
\end{equation*}
and so 
\begin{equation*}
\hom(T,R')=h=\hom(S,R').
\end{equation*}
We restrict to the open subvariety of $R(Q,h(a\ b))$ where
\begin{equation*}
\hom(T,R_{p})=h=\hom(S,R_{p}).
\end{equation*}
Since the homomorphism from
$\Hom(T,R^{h})\otimes T$ to $R^{h}$ is surjective there is an open
subvariety where the natural homomorphism from $\Hom(T,R_{p})\otimes
T$ to $R_{p}$ is surjective. Let $W_{p}$ be the kernel of the natural
homomorphism from $\Hom(S,T)\otimes \Hom(T,R_{p})$ to $\Hom(S,R_{p})$
which has dimension $hp$ on this open subvariety. Since the natural
homomorphism from $W_{p}\otimes S$ to $\Hom(T,R_{p})\otimes T$ has
image in the kernel of the homomorphism to $R_{p}$ and is an
isomorphism with the kernel when $R_{p}\cong R^{h}$ it follows that on
a suitable open subvariety the complex
\begin{equation*}
\ses{W_{p}\otimes S}{\Hom(T,R_{p})\otimes T}{R_{p}}
\end{equation*}
is a short exact sequence which proves the third part of the lemma.
\end{proof}

\begin{theorem}
\label{th:fundreg}
A Schur root for the $n$th Kronecker quiver $Q$ in the fundamental
region for the action of the Weyl group is reducible to matrix normal
form.
\end{theorem}
\begin{proof}
Let the Schur root be $h(a\ b)$ where $(a\ b)$ is indivisible and
$a\leq b$ which we may assume by lemma
\ref{lem:dual}.  Let $S$ and $T$ be the representations constructed in
the preceding paragraphs and let $p=1-\langle(a\ b),(a\ b)
\rangle$. Let $Q'$ be the $(1+p)$th Kronecker quiver. Then we show
that for a general representation $R'$ of dimension vector $h(a\ b)$,
\begin{equation*}
R'\cong \Hom_{Q}(S\oplus T,R')\otimes_{Q'}(S\oplus T)
\end{equation*}
and for a
general representation $R''$ of dimension vector $(h,h)$ for $Q'$,
\begin{equation*}
R''\cong \Hom_{Q}(S\oplus T,R''\otimes_{Q'}(S\oplus T)).
\end{equation*}

Note that if $P_{1}$ and $P_{2}$ are the projective representations of
dimension vector $(0\ 1)$ and $(1\ p+1)$ for the quiver $Q'$ then
$P_{1}\otimes_{Q'}(S\oplus T)\cong S$ whilst
$P_{2}\otimes_{Q'}(S\oplus T)\cong T$.

Let $R'$ be a general representation of dimension vector $h(a\ b)$ for
the $n$th Kronecker quiver. Then we know that
$\hom(S,R')=h=\hom(T,R')$ and the kernel of the natural homomorphism
from $\Hom(T,R')\otimes T$ to $R'$ which is surjective is isomorphic
to $S^{hp}$. In fact, let $W'$ be the kernel of the homomorphism from
$\Hom(S,T)\otimes \Hom(T,R')$ to $\Hom(S,R')$. Then for general $R'$
we showed in the last lemma that 
\begin{equation*}
\ses{W'\otimes S}{\Hom(T,R')\otimes T}{R'}
\end{equation*}
is a short exact sequence. However, 
\begin{equation*}
\ses{W'\otimes P_{1}}{\Hom(T,R')\otimes P_{2}}{\Hom_{Q}(S\oplus T,R')}
\end{equation*}
is the projective resolution of $\Hom_{Q}(S\oplus T,R')$; tensoring
the second of these short exact sequences by $S\oplus T$ shows that
\begin{equation*}
R'\cong \Hom_{Q}(S\oplus T,R')\otimes_{Q'}(S\oplus T).
\end{equation*}

Conversely, let $W$ be a general vector subspace of dimension $hp$ in
$k^{h}\otimes \Hom(S,T)$ so that 
\begin{equation*}
\ses{W\otimes P_{1}}{k^{h}\otimes P_{2}}{R}
\end{equation*}
is the projective resolution of a general representation
$R$ of dimension vector $(h\ h)$ for the quiver $Q'$. Then since there
exist such subspaces for which the homomorphism from $W\otimes S$ to
$k^{h}\otimes T$ is injective, this remains true for a general $W$ and
the cokernel of this homomorphism is a representation $R_{W}$ of
dimension vector $h(a\ b)$ for the quiver $Q$. Since there exists a
choice of $W$ for which $\hom(S,R_{W})=h=\hom(T,R_{W})$, this remains
true for a general $W$, and so 
\begin{equation*}
\ses{W\otimes P_{1}}{k^{h}\otimes P_{2}}{\Hom(S\oplus T,R_{W})}
\end{equation*}
is the projective resolution of $\Hom(S\oplus T,R_{W})$ (noting that
$\Hom(T,R_{W})$ is naturally isomorphic to $k^{h}$) which means that
$R$ is isomorphic as required to $\Hom_{Q}(S\oplus
T,R\otimes_{Q'}(S\oplus T))$.
\end{proof}

\begin{theorem}
\label{th:KQreduction}
Every Schur root for the $n$th Kronecker quiver is reducible to matrix
normal form.
\end{theorem}
\begin{proof}
This follows at once from the preceding theorem, theorem \ref{KQ1} and
lemma \ref{lem:reflection}. 
\end{proof}

\section{Internal structure of Schur representations}
\label{internal}

The main aim of this section is to find information on the internal
structure of a general representation of dimension vector $\alpha$. In
the case where $\alpha$ is a Schur root it allows us to show that
$\alpha$ is built up from two smaller Schur roots in a useful way and
for a general dimension vector it takes the form of an algorithm that
computes the canonical decomposition. 
 
For the moment we make the assumption that $Q$ is a quiver without
loops. Firstly recall some facts about the canonical decomposition. The sum
$\alpha=\sum_{i}\beta_{i}$ is the canonical decomposition if and only
if each $\beta_{i}$ is a Schur root and $\ext(\beta_{i},\beta_{j})=0$
if $i\neq j$. If $\beta_{i}=\beta_{j}$ and $i\neq j$ then either
$\beta_{i}$ is a real Schur root, that is,
$\hom(\beta_{i},\beta_{i})=1$, or else $\hom(\beta_{i},\beta_{i})=0$
that is, $\beta_{i}$ is an \definition{isotropic} Schur root. If
$\beta_{i}\neq \beta_{j}$ then one of $\hom(\beta_{i},\beta_{j})$ and
$\hom(\beta_{j},\beta_{i})$ must be $0$. In fact the following lemma
contains a stronger statement.

\begin{lemma}
\label{lem:noloops}
Let $\{\beta_{i}:i=1\rightarrow m\}$ be pairwise distinct canonical
summands of the dimension vector $\alpha$. Then one of
$\hom(\beta_{i},\beta_{i+1})$ for $i=1\rightarrow m-1$ and
$\hom(\beta_{m},\beta_{1})$ must be $0$. 
\end{lemma}
\begin{proof}
Let $\{R_{i}:i=1\rightarrow m\}$ be representations of dimension
vector $\beta_{i}$ respectively such that $\Ext(R_{i},R_{j})=0$ for
$i\neq j$. Then by lemma 4.1 of \cite{hapring}, any homomorphism
from $R_{i}$ to $R_{j}$ must be either injective or surjective. Assume
our conclusion is false; then there is a non-zero homomorphism from
$R_{i}$ to $R_{i+1}$ for each $i<m$ and a non-zero homomorphism from
$R_{m}$ to $R_{1}$. Each homomorphism in this chain must be either
surjective or injective but not both and no surjective homomorphism
may be followed by an injective homomorphism since their composition
would then be neither injective nor surjective; hence these
homomorphisms must be either all surjective or all injective which is
absurd.
\end{proof}

This allows us to prove the following useful lemma.

\begin{lemma}
\label{lem:min+and}
Let $\alpha=\sum_{i}n_{i}\beta_{i}$ be the canonical decomposition of
$\alpha$ where $\beta_{i}=\beta_{j}$ if and only if $i=j$. Then it is
possible to choose the indexing so that $i<j\Rightarrow
\hom(\beta_{i},\beta_{j})=0$.  
\end{lemma}
\begin{proof}
We define a relation $<$ on the dimension vectors $\{\beta_{i}\}$ in
the canonical decomposition by $i<j$ if
$\hom(\beta_{j},\beta_{i})>0$. Then the preceding lemma shows that $<$
is a partial order and can therefore be extended to a total order
which is the conclusion of the lemma.
\end{proof}

Thus one can assume that the canonical decomposition of $\alpha$ is 
$\alpha=\sum_{i}n_{i}\beta_{i}$ where $\beta_{i}=\beta_{j}$ if and
only if $i=j$ and $\hom(\beta_{1},\beta_{i})=0$ for $i>1$. 

If $\alpha$ is an indivisible Schur root then the canonical
decomposition of $\beta=n\alpha$ is $\beta$ if $\alpha$ is neither
real nor isotropic and in these two cases it is $n\alpha$ as shown in
theorem 3.8 of \cite{genrep}. If $\beta=n\alpha$ where $\alpha$ is an
indivisible Schur root we say that $\beta$ is a \definition{uniform}
dimension vector and that $\alpha$ is the \definition{root} of $\beta$
when $\alpha$ is real or isotropic and that $\beta$ is its own root if
$\alpha$ (and hence $\beta$) is non-isotropic; we shall use the
notation $(\beta)$ for the canonical decomposition of such a uniform
dimension vector. Thus we can re-write the canonical decomposition of
a general dimension vector in the form $\alpha=\sum_{i}(\beta_{i})$
where each $\beta_{i}$ is a uniform dimension vector and the root of
$\beta_{i}$ does not equal the root of $\beta_{j}$ when $i\neq
j$. Then by theorem 3.8 of \cite{genrep} it follows that the canonical
decomposition of $n\alpha$ is $n\alpha=\sum_{i}(n\beta_{i})$.

In \cite{genrep}, the author showed in characteristic $0$ that every
representation of dimension vector $\alpha$ contains a
subrepresentation of dimension vector $\beta$ if and only if
$\ext(\beta,\alpha-\beta)=0$ and this was extended to arbitrary
characteristic in \cite{Craw}. We shall use this result repeatedly in
this section.  A dimension vector $\beta$ is said to be a
\definition{rigid sub-dimension vector} of the dimension vector
$\alpha$ if and only if a general representation of dimension vector
$\alpha$ has a unique subrepresentation of dimension vector
$\beta$. We shall be interested in finding uniform rigid sub-dimension
vectors of a particular dimension vector $\alpha$; we may as well
assume that the support of $\alpha$ is $Q$. When $Q$ is a quiver
without loops then there are obviously such dimension vectors since if
$v$ is a sink vertex, then $\alpha(v)e_{v}$ is a uniform rigid
sub-dimension vector of $\alpha$ where $e_{v}$ is the dimension vector
such that $e_{v}(v)=1$ and $e_{v}(w)=0$ when $v\neq w$.

\begin{lemma}
\label{lem:unirig+and}
Let $\alpha=\sum_{i}(\beta_{i})$ be the canonical decomposition of
$\alpha$ where each $\beta_{i}$ is a uniform dimension vector such
that the root of each $\beta_{i}$ is distinct. Assume that
$\hom(\beta_{1},\beta_{i})=0$ for $i\neq 1$. Then the dimension vector
$\beta_{1}$ is a uniform rigid sub-dimension vector of $\alpha$.
\end{lemma}
\begin{proof}
To begin with, note that
$\hom(\beta_{1},\beta_{i})=0=\ext(\beta_{1},\beta_{i})$ for $i>1$ and
hence $\hom(\beta_{1},\alpha-\beta_{1})=0
=\ext(\beta_{1},\alpha-\beta_{1})$. Recall from section 3 of
\cite{genrep} the algebraic variety $R(Q,\beta_{1}\subset
\alpha)$ which parametrises representations of dimension vector
$\alpha$ with a distinguished subrepresentation of dimension vector
$\beta_{1}$. Then from section 3 of \cite{genrep} it follows that
the morphism from $R(Q,\beta_{1}\subset \alpha)$ to $R(Q,\alpha)$
is surjective since the fibre above a point $p$ is bijective with the
set of subrepresentations of dimension vector $\beta_{1}$ in
$R_{p}$. Moreover, the conditions above mean that the dimension of
$R(Q,\beta_{1}\subset \alpha)$ equals the dimension of
$R(Q,\alpha)$ so that the fibre above a general point is finite. Both
varieties are irreducible. Moreover, there is a rational section, a
morphism defined on an open subvariety of $R(Q,\alpha)$ to
$R(Q,\beta_{1}\subset \alpha)$ which sends the point $p$ to the
point in the fibre above $p$ corresponding to the subrepresentation
that is the direct summand of dimension vector $\beta_{1}$. But a
morphism between irreducible algebraic varieties that is generically
finite and has a rational section must be generically bijective
(consider the effect of the rational section on the function fields of
the two varieties). This implies that the fibre above a general point
of $R(Q,\alpha)$ consists of $1$ point which means that
$\beta_{1}$ is a rigid sub-dimension vector of $\alpha$ and it is
uniform by assumption.
\end{proof}

Let $\alpha=\sum_{i}(\gamma_{i})$ be the canonical decomposition of
$\alpha$ where each $\gamma_{i}$ is a uniform dimension vector and the
root of $\gamma_{i}$ equals the root of $\gamma_{j}$ if and only if
$i=j$. If $j$ is an index such that $\hom(\gamma_{j},\gamma_{i})=0$
for $i\neq j$ then call the dimension vector $\gamma_{j}$ a
\definition{uniform rigid summand} of $\alpha$. Since the canonical
decomposition of $n\alpha$ is $n\alpha=\sum_{i}(n\gamma_{i})$ it
follows that $\gamma_{j}$ is a uniform rigid summand of $\alpha$ if
and only if $n\gamma_{j}$ is a uniform rigid summand of
$n\alpha$. Since we shall make use of this fact later we note it in
the following form.

\begin{lemma}
\label{lem:sumhcf}
Let $\alpha$ be a dimension vector for the quiver $Q$. Then
$g(\alpha)=\hcf_{v}(\alpha(v))$ divides $g(\gamma)$ if $\gamma$ is a
uniform rigid summand of $\alpha$ 
\end{lemma}

The next two lemmas give ways to find new rigid sub-dimension vectors
of a dimension vector which will form the basis of an algorithm to
compute the canonical decomposition of a dimension vector $\alpha$ or
else to find a rigid sub-dimension vector $\beta$ of $\alpha$ when it
is a Schur root such that both $\beta$ and $\alpha-\beta$ are uniform.

\begin{lemma}
\label{lem:swap1}
Let $\alpha=\beta+\gamma+\delta$ where $\beta+\gamma$ is a rigid
sub-dimension vector of $\alpha$ and $\beta$ is a uniform rigid
summand of $\beta+\gamma$. Then $\beta$ is a rigid sub-dimension
vector of $\alpha$.
\end{lemma}
\begin{proof}
Let $R$ be a general representation of dimension vector $\alpha$; in
particular, it has a unique subrepresentation $R'$ of dimension vector
$\beta+\gamma$ which in turn may be taken to have a unique
subrepresentation $R''$ of dimension vector $\beta$. Since $\beta$ is
a uniform rigid summand of $\beta+\gamma$, $\ext(\beta,\gamma)=0$ and
also $\ext(\beta,\delta)=0$ since $\ext(\beta+\gamma,\delta)=0$ and
$\beta$ is a canonical summand of $\beta+\gamma$. So,
$\ext(\beta,\gamma+\delta)=0$ and hence every representation of
dimension vector $\alpha$ has a subrepresentation of dimension vector
$\beta$. Also since $\gamma$ is a canonical summand of $\beta+\gamma$
and $\ext(\beta+\gamma,\delta)=0$ then $\ext(\gamma,\delta)=0$ so a
general representation of dimension vector $\gamma+\delta$ has a
subrepresentation of dimension vector $\gamma$. So let $S$ be a
subrepresentation of $R$ of dimension vector $\beta$; then $R/S$ has a
subrepresentation of dimension vector $\gamma$, $T/S$ where $T$ has to
be $R'$ since it is a subrepresentation of dimension vector
$\beta+\gamma$ of $R$. So, $S$ is a subrepresentation of dimension
vector $\beta$ in $R'$ and must be $R''$. Thus $\beta$ is a rigid
sub-dimension vector of $\alpha$.
\end{proof}
\begin{lemma}
\label{lem:swap2}
Let $\alpha=\beta+\gamma+\delta$ where $\beta$ is a rigid
sub-dimension vector of $\alpha$ and $\gamma$ is a uniform rigid
summand of $\gamma+\delta$. Then $\beta+\gamma$ is a rigid
sub-dimension vector of $\alpha$.
\end{lemma}
\begin{proof}
Since $\ext(\beta,\gamma+\delta)=0$ and $\gamma$ is a uniform rigid
summand of $\gamma+\delta$,
$\ext(\beta,\gamma)=0=\ext(\beta,\delta)$. Also
$\ext(\gamma,\delta)=0$ and so $\ext(\beta+\gamma,\delta)=0$. Hence a
general representation of dimension vector $\alpha$ has a
subrepresentation of dimension vector $\beta+\gamma$. In turn this has
a subrepresentation of dimension vector $\beta$ which must be the
unique one and the factor must be the unique subrepresentation of
dimension vector $\gamma$. So there is a unique subrepresentation of
dimension vector $\beta+\gamma$.
\end{proof}
These results are the basis for the following lemma which contains
most of the work for our understanding of the internal structure of
Schur roots and representations and is also the basis for an algorithm
for constructing the canonical decomposition of a dimension vector.
\begin{lemma}
\label{t5}
Let $\alpha=\beta+\gamma+\delta$ be a dimension vector for the quiver
$Q$ without loops where $\beta$ is a uniform rigid
sub-dimension vector of $\alpha$ and $\gamma$ is a uniform rigid
summand of $\gamma+\delta$. Let $\beta=m\beta'$ and $\gamma=n\gamma'$
be the canonical decompositions of $\beta$ and $\gamma$ so that
$\beta'$ is the root of $\beta$. Then either $\gamma'$ is a canonical
summand of $\alpha$ or else any uniform rigid summand of
$\beta+\gamma$ is a uniform rigid sub-dimension vector of $\alpha$
whose root is larger than $\beta'$, the root of $\beta$.
\end{lemma}
\begin{proof}
Our assumptions imply that
\begin{equation*}
\hom(\beta,\gamma+\delta)=0=\ext(\beta,\gamma+\delta)
\end{equation*}
and hence
\begin{align*}
\hom(\beta,\gamma)=&0=\ext(\beta,\gamma)\\
\hom(\beta,\delta)=&0=\ext(\beta,\delta)\\
\ext(\gamma,\delta)=&0=\ext(\delta,\gamma)
\end{align*}
since $\gamma$ and $\delta$
are canonical summands of $\gamma+\delta$. If
$\ext(\gamma,\beta)=0$ then $\gamma$ is a summand of $\alpha$ because
$\alpha-\gamma=\beta+\delta$ and so
\begin{equation*}
\ext(\gamma,\alpha-\gamma)=0=\ext(\alpha-\gamma,\gamma).
\end{equation*}
Thus if $\gamma'$ is not a canonical summand of $\alpha$ then
$\ext(\gamma,\beta)>0$ and therefore $\ext(\gamma',\beta')>0$. 
It follows that
$\hom(\gamma',\beta')=0$ since $\ext(\beta',\gamma')=0$ and so
one of $\hom(\gamma',\beta')$ and $\ext(\gamma',\beta')$ must be
$0$ by theorem 4.1 of \cite{genrep}.

By lemma \ref{lem:swap2}, $\beta+\gamma$ is a rigid sub-dimension
vector of $\alpha$.  So let us consider its canonical
decomposition. Firstly 
\begin{equation*}
\hom(\beta',\gamma')=0=\ext(\beta',\gamma')
\end{equation*}
and $\hom(\gamma',\beta')=0$. On the other hand,
$\ext(\gamma',\beta')\neq 0$. Let $S$ be a general representation of
dimension vector $\beta+\gamma$. Then $S$ has a subrepresentation $R$
of dimension vector $\beta$ such that both $R$ and $S/R$ are general
representations. So $R\cong \oplus_{j=1}^{m}R_{j}$ and $S/R\cong
\oplus_{i=1}^{n}S_{i}$ where $\ddim R_{j}=\beta'$, $\ddim
S_{i}=\gamma'$, each $R_{j}$ and $S_{i}$ is a Schur representation,
\begin{equation*}
\Hom(R_{j},S_{i})=0=\Hom(S_{i},R_{j})
\end{equation*}
and also $\Hom(R_{j},R_{l})=0$ for
$j\neq l$ unless $\beta'$ is a real Schur root in which case they are
isomorphic, and similarly, $\Hom(S_{i},S_{l})=0$ for $i\neq l$ unless
$\gamma'$ is a real Schur root in which case they are isomorphic. By
Ringel's simplification process \cite{ringel} one knows that any
summand of $S$ must have a filtration by subrepresentations such that
the factors are isomorphic to either an $S_{i}$ or an $R_{j}$. No
$S_{i}$ can be a summand of $S$ since then its dimension vector
$\gamma'$ would be a summand in the canonical decomposition of
$\alpha$. Thus any summand in the canonical decomposition of
$\beta+\gamma$ must be of the form $a\beta'+b\gamma'$ where $a>0$. If
$\beta'$ is not a summand in the canonical decomposition of
$\beta+\gamma$, then every canonical summand of $\beta+\gamma$ is
larger than $\beta'$.  If $\beta'$ is a summand in the canonical
decomposition of $\beta+\gamma$, we shall soon see that it cannot be
the root of a uniform rigid summand. Thus it follows that $\langle
(m-1)\beta'+n\gamma',\beta'\rangle\geq 0$ and so $\langle
\beta',\beta'\rangle>0$ which implies that $\beta'$ is a real Schur
root and so each $R_{j}$ is isomorphic to the real Schur
representation $G(\beta')$ of dimension vector $\beta'$. However,
there is a summand of $S$ that has a proper subrepresentation
isomorphic to some $R_{j}\cong G(\beta')$ and so, as claimed, $\beta'$
cannot be a uniform rigid summand of $S$. It follows that a uniform
rigid summand of $\beta+\gamma$ cannot be a multiple of $\beta'$ in
this case and therefore the root of a uniform rigid summand $\kappa$
of $\beta+\gamma$ must be larger than $\beta'$.  Therefore, in either
case there is a uniform rigid summand $\kappa$ of $\beta+\gamma$ whose
root is larger than $\beta'$.  However, by lemma \ref{lem:swap1},
$\kappa$ is a uniform rigid sub-dimension vector of $\alpha$ which
completes our proof.
\end{proof}
The main result on the internal structure of Schur roots is now a
simple induction.
\begin{theorem}
\label{t6}
Let $\alpha$ be a Schur root for the quiver $Q$ without loops. Then
there exists a rigid sub-dimension vector $\beta$ of $\alpha$ such
that both $\beta$ and $\alpha-\beta$ are uniform dimension vectors. 
More particularly, if $\epsilon$ is a uniform rigid sub-dimension
vector of $\alpha$ such that $g(\alpha)|g(\epsilon)$ then we may
choose $\beta$, a uniform rigid sub-dimension dimension vector of
$\alpha$ such that $\alpha-\beta$ is uniform, the root of $\beta$ is
at least as large as that of $\epsilon$ and
$\hcf(g(\beta),g(\alpha-\beta))=g(\alpha)$. 
\end{theorem}
\begin{proof}
We may assume that the support of $\alpha$ is $Q$. Let $v$ be a sink vertex;
then as discussed before $\alpha(v)e_{v}$ is a uniform rigid
sub-dimension vector of $\alpha$ where $e_{v}$ is the dimension vector
of the simple representation at the vertex $v$. Therefore, if we do
not already have a dimension vector $\epsilon$ we may take
$\epsilon=\alpha(v)e_{v}$. 

If the factor is uniform we simply need to check the numerical
statements. But $g(\alpha)|g(\epsilon)$ implies that
$g(\alpha)|\hcf(g(\epsilon),g(\alpha-\epsilon))|g(\alpha)$ and so
equality follows.

Otherwise, assume that $\epsilon$ is a uniform rigid sub-dimension
vector of $\alpha$ such that $\alpha-\epsilon$ is not uniform and
proceed by induction on $\alpha-\epsilon'$ where $\epsilon'$ is the
root of $\epsilon$. Since $\alpha-\epsilon$ is not uniform,
$\alpha-\epsilon=\gamma+\delta$ where $\gamma$ is a uniform rigid
summand of $\alpha-\epsilon$. Since the root of $\gamma$ cannot be a
summand of $\alpha$ one concludes by lemma \ref{t5} that there is a
uniform rigid sub-dimension vector $\epsilon_{1}$ of $\alpha$ whose
root is larger than $\epsilon'$. Since $g(\alpha)|g(\alpha-\epsilon)$
it follows that $g(\alpha)|g(\gamma)$ by lemma \ref{lem:sumhcf} hence
$g(\alpha)|g(\epsilon+\gamma)$ and hence $g(\alpha)$ must divide
$\epsilon_{1}$ by lemma \ref{lem:sumhcf} since it is a uniform rigid
summand of $\epsilon+\gamma$.  By induction, the result follows.
\end{proof}

Lemma \ref{t5} also gives an algorithm to compute the canonical
decomposition of a dimension vector $\alpha$ in the following
way. First of all, it is a simple matter to compute whether a
dimension vector is a Schur root on a quiver with $2$ vertices $v$ and
$w$ with no arrows from $w$ to $v$ though all other possible arrows
are allowed, since all roots are Schur. Now assume that $Q$ is a
quiver without loops.  As at the beginning of the proof of theorem
\ref{t6}, one may assume that the support of $\alpha$ is the quiver
$Q$ and $v$ is a sink vertex so that $\alpha(v)e_{v}$ is a uniform
rigid sub-dimension vector of $\alpha$ and $\alpha-\alpha(v)e_{v}$ is
a smaller dimension vector. One can compute the canonical
decomposition of $\alpha-\alpha(v)e_{v}$ by induction. If
$\alpha-\alpha(v)e_{v}$ is not uniform then lemma \ref{t5} either
gives a canonical summand of $\alpha$ and we may proceed by induction
or else it gives a new uniform rigid sub-dimension vector of $\alpha$
which has a larger root. Thus the only time a problem occurs is when
there is a uniform rigid sub-dimension vector $\beta$ of $\alpha$ such
that $\gamma=\alpha-\beta$ is also uniform. Let $\beta=m\beta'$ and
$\gamma=n\gamma'$ be their canonical decompositions.  If both $\beta'$
and $\gamma'$ are not real then it is a simple matter to see that
$\alpha$ must be a Schur root and $\beta$ is a uniform rigid
sub-dimension vector of $\alpha$ such that $\alpha-\beta$ is also
uniform. Without loss of generality one may assume that $\beta'$ is
real (for example by reversing all the arrows of the quiver). Let
$t=\ext(\gamma',\beta')$. If $\gamma'$ is also real, one considers the
canonical decomposition of $(n,m)$ for the $2$-vertex quiver with $t$
arrows from the first vertex to the second and no loops. If this is
$(n,m)=a(b,c)+d(e,f)$ then the canonical decomposition of $\alpha$ is
$\alpha=a(b\gamma+c\beta)+d(e\gamma+f\beta)$. Otherwise, $\alpha$ is a
Schur root. If $\gamma'$ is isotropic, one again considers the
dimension vector $(n,m)$ for the $2$-vertex quiver with $t$ arrows
from the first to the second vertex but with $1$ loop at the first
vertex. Again its canonical decomposition determines the canonical
decomposition of $\alpha$ using the same formula. Finally, if
$\gamma'$ is neither real nor isotropic then $\gamma=\gamma'$ and
$\alpha$ is a Schur root if and only if $m\leq t$; otherwise its
canonical decomposition is $\alpha=(\gamma+t\beta')+(m-t)\beta'$.
Note that when $\alpha$ is a Schur root we have also calculated a
dimension vector $\beta$ that is a rigid sub-dimension vector of
$\alpha$ such that both $\beta$ and $\alpha-\beta$ are uniform.

In \cite{genrep}, there is a different algorithm for computing the
canonical decomposition of a dimension vector. One knows that $\beta$
is a canonical summand of the dimension vector $\alpha$ if and only if
it is a Schur root and $\ext(\beta,\alpha-\beta)=0
=\ext(\alpha-\beta,\beta)$ or equivalently a representation of
dimension vector $\alpha$ has subrepresentations of dimension vector
$\beta$ and $\alpha-\beta$. On the other hand,
$\ext(\beta,\alpha-\beta)$ and $\ext(\alpha-\beta,\beta)$ may be
calculated by knowing the dimension vectors of subrepresentations of
$\beta$. Since $\beta<\alpha$ inductively we know all of these
dimension vectors. This algorithm is substantially more complicated
than the one in this paper and is correspondingly much slower.

One should also note that this extends to computing the canonical
decomposition of a dimension vector over an arbitrary quiver. If $Q$
is a quiver with vertex set $V$ and arrow set $A$ and $\alpha$ is a
dimension vector for the quiver $Q$ then one constructs a new quiver,
the \definition{double} of $Q$, $Q'$, with vertex set $V\times
\{0,1\}$, and arrow set $V\times\{2\}\cup (A\times \{0\})$ where
$i(v,2)=(v,0)$, $t(v,2)=(v,1)$, $i(a,0)=(ia,0)$ and
$t(a,0)=(ta,1)$. Then $\alpha'=\alpha p_{1}$, where $p_{1}$ is
projection on the first factor $V$ of $V\times \{0,1\}$, is a
dimension vector on $Q'$ and if $\alpha=\sum_{i}\beta_{i}$ is the
canonical decomposition of $\alpha$ then the canonical decomposition
of $\alpha p_{1}$ is $\alpha'=\sum_{i}\beta_{i}'$ since for a general
representation of dimension vector $\alpha'$, the arrows in $V$ are
invertible and the full subcategory of representations of $Q'$ such
that these arrows are invertible is naturally equivalent to the
category of representations of the quiver $Q$. To see this quickly we
note that given a representation $R$ of $Q$, we associate a
representation $R'$ of $Q'$ by defining $R'(v,i)=R(v)$ for $i=0,1$,
$R'(v,2)$ is the identity from $R(v)$ to itself and
$R'(a,0)=R(a)$. Conversely given a representation $S$ of the quiver
$Q'$ such that each $R(v,2)$ is invertible, we define a representation
$S'$ of the quiver $Q$ by $S'(v)=S(v,0)$ and
$S'(a)=S(a,o)S(v,2)^{-1}$. It is simple to check that these
assignments define functors to demonstrate the natural equivalence
claimed.

\section{Moduli spaces of representations of quivers}
\label{finally}

Let $Q(p,t,q)$ be the quiver with two vertices $v$ and $w$, $p$ loops
at the vertex $v$, $q$ loops at the vertex $w$ and $t$ arrows from $v$
to $w$. The first theorem of this section shows roughly speaking that
if every Schur root for the quiver $Q(p,t,q)$ is reducible to matrix
normal form then every Schur root for a quiver without loops is also
reducible to matrix normal form. The rest of the section proves this
is true for these special two vertex quivers and theorem
\ref{th:complete} of this section uses the double quiver introduced at
the end of the last section to show that every Schur root for every
quiver is reducible to matrix normal form. The last theorem spells out
those moduli spaces of representations we are now able to show to be
rational varieties.

\begin{theorem}
\label{th:reduce}
Let $\alpha$ be a Schur root for the quiver $Q$ and let $n\gamma$ be a
uniform rigid sub-dimension vector such that $m\beta=\alpha-n\gamma$
is also uniform where the roots of $m\beta$ and $n\gamma$ are Schur
roots that are reducible to matrix normal form and $\beta$ and
$\gamma$ are indivisible. Then if the dimension vector $(m\ n)$ for
the quiver $Q(p,t,q)$, where $p=p(\beta)$, $q=p(\gamma)$ and
$t=\ext(\beta,\gamma)$, is reducible to matrix normal form, so is
$\alpha$.
\end{theorem}
\begin{proof}
Let $PGl_{m,n}$ be the factor of $Gl_{m}\times Gl_{n}$ by the diagonal
embedding of $k^{*}$. The first step of the proof is to construct a
$(Gl_{m}\times Gl_{n},PGl_{m,n})$-standard family of representations
of dimension vector $\alpha$ for the quiver $Q$ which is $PGl_{m,n}$
birational to $R(Q(p,t,q),(m\ n))$. 

The first thing we need to do this are suitable familes of
representations of dimension vector $m\beta$ and $n\gamma$. We deal
only with $m\beta$ since the same construction is needed for
$n\gamma$. We need a family $\mathcal{R}_{m\beta}$ over the algebraic
variety $X_{m\beta}$ such that $Gl_{m}$ acts on $\mathcal{R}_{m\beta}$
compatibly with an action of $PGl_{m}$ on $X_{m\beta}$ and the
stabiliser in $Gl_{m}$ of a point in $X_{m\beta}$ acts on the
corresponding representation as the group of automorphisms of that
representation. 

If $\beta$ is a real Schur root then there is a unique real
Schur representation $G(\beta)$ of dimension vector $\beta$. We take
the variety $X_{m\beta}$ to be a point on which $PGl_{m}$ acts
trivially and the family $\mathcal{R}_{m\beta}$ will be simply the
representation $k^{m}\otimes G(\beta)$ on which $Gl_{m}$ acts via its
action on $k^{m}$. If $\beta$ is non-isotropic then the root of
$m\beta$ is itself and our assumptions imply that that there is a
$(Gl_{m},PGl_{m})$-standard family $\mathcal{R}_{m\beta}$ of
representations of dimension vector $m\beta$ over an algebraic variety
$X_{m\beta}$ which is $PGl_{m}$ birational to $M_{m}(k)^{p}$. 

In the case where $\beta$ is isotropic there is more to do. In this
case, our assumptions imply that there is an algebraic variety
$X_{\beta}$ which is an open subvariety of $\mathbb{A}^{1}$, the
affine line, which carries a family $\mathcal{R}_{\beta}$ of
representations of dimension vector $\beta$ and that this family is
$(k^{*},\{1\})$-standard which in this situation simply means that
different points of $X_{\beta}$ give non-isomorphic
representations. We wish to construct from this a family of
representations of dimension vector $m\beta$ on which $Gl_{m}$ acts,
defined over an algebraic variety $X_{m\beta}$ on which $PGl_{m}$ acts
so that $X_{m\beta}$ is $PGl_{m}$ birational to an open subvariety of
$M_{m}(k)$. We take $Y_{m\beta}$ to be the open subvariety of
$X_{\beta}^{m}$ consisting of $m$ distinct ordered points of
$X_{\beta}$. This carries in the obvious way a family
$\mathcal{S}_{m\beta}$ of representations of dimension vector
$m\beta$. The natural action of the symmetric group $S_{m}$ on
$Y_{m\beta}$ has the property that its orbits correspond to the
isomorphism classes of representations in the family
$\mathcal{S}_{m\beta}$. Of course, this action of $S_{m}$ on
$Y_{m\beta}$ lifts to an action on $\mathcal{S}_{m\beta}$ but this is
not quite all we need; the algebraic group ${k^{*}}^{m}$ acts on
$\mathcal{S}_{m\beta}$ via the action of $k^{*}$ on
$\mathcal{R}_{\beta}$ and the algebraic group $H$ generated by
${k^{*}}^{m}$ and $S_{m}$ acts on $\mathcal{S}_{m\beta}$. The subgroup
${k^{*}}^{m}$ of $H$ is normal with factor group $S_{m}$ and we shall
regard $H$ as acting on $Y_{m\beta}$ so that ${k^{*}}^{m}$ acts
trivially and $S_{m}$ acts as described. Thus the actions of $H$ on
$\mathcal{S}_{m\beta}$ and $Y_{m\beta}$ are compatible. Now $H$ is a
subgroup of $Gl_{m}$ in the obvious way; it is the stabiliser with
respect to the action of $Gl_{n}$ by conjugation of the vector
subspace of $M_{m}(k)$ consisting of diagonal matrices. Therefore, we
can form the family
$\mathcal{R}_{m\beta}=\mathcal{S}_{m\beta}\times^{H}Gl_{m}$ of
representations of dimension vector $m\beta$ over the algebraic
variety $X_{m\beta}=Y_{m\beta}\times^{H}Gl_{m}$. An identification of
$Y_{m\beta}$ with some open subvariety of the vector space of diagonal
matrices shows that $X_{m\beta}$ is $PGl_{m}$ birational to
$M_{m}(k)$, each $\mathcal{R}_{m\beta}(v)$ is a $Gl_{m}$ vector bundle
of weight $1$ and for each point $p$ of $X_{m\beta}$, the stabiliser
in $Gl_{m}$ of $p$ acts on the corresponding representation as the
group of automorphisms of that representation since that is true for
points in $Y_{m\beta}$.

Note that there is a short exact sequence of groups
\begin{equation*}
1\rightarrow k^{*}\rightarrow PGl_{m,n}\rightarrow PGl_{m}\times
PGl_{n}\rightarrow 1.
\end{equation*}
Thus a variety on which $PGl_{m}\times PGl_{n}$ acts may
have $PGl_{m,n}$ vector bundles of weight $1$.

Let $\mathcal{R}=\mathcal{R}_{m\beta}$ and
$\mathcal{S}=\mathcal{R}_{n\gamma}$ be families of representations of
dimension vector $m\beta$ and $n\gamma$ over the algebraic varieties
$X_{m\beta}$ and $X_{n\gamma}$ as constructed above.  So $X_{m\beta}$
is $PGl_{m}$-birational to $M_{m}(k)^{p}$ and $X_{n\gamma}$ is
$PGl_{n}$-birational to $M_{n}(k)^{q}$. There exists a non-empty open
$PGl_{m}\times PGl_{n}$-equivariant subvariety $U$ of $X_{m\beta}\times
X_{n\gamma}$ consisting of the points $(x,y)$ where
$\ext(\mathcal{R}_{x},\mathcal{S}_{y})=mnt$.  Then over $U$ there is a
$PGl_{m,n}$ vector bundle $E$ of weight $1$ whose fibre above the
point $(x,y)$ is $\Ext(\mathcal{R}_{x},\mathcal{S}_{y})$. However,
$R(Q(p,t,q),(m,n))$ is a $PGl_{m,n}$ vector bundle of weight $1$ over
$M_{m}(k)^{p}\times M_{n}(k)^{q}$ and by lemma \ref{lem:tg} it follows
that $E$ is $PGl_{m,n}$ birational to $R(Q(p,t,q),(m,n))$.

Further $E$ carries a family $\mathcal{T}$ of representations of
dimension vector $\alpha$ of the quiver $Q$ and it is clear that
$\mathcal{T}$ is a general family since every representation of
dimension vector $\alpha$ has a subrepresentation of dimension vector
$n\gamma$ and it is an open condition that one such subrepresentation
should be isomorphic to $\mathcal{S}_{y}$ for some point $y\in
X_{n\gamma}$ whilst the factor should be isomorphic to
$\mathcal{R}_{x}$ for some point $x\in X_{m\beta}$. Since a general
representation of dimension vector $\alpha$ has a unique
subrepresentation of dimension vector $n\gamma$, it follows that
$E'=\{z\in E:\mathcal{T}_{z} \text{ has a unique subrepresentation of
dimension vector }n\gamma\}$ is a non-empty open subvariety of $E$ and
we shall see that two points of $E'$ give rise to isomorphic
representations if and only if they are in the same orbit for the
action of $PGl_{m,n}$. To show this we represent the points of $E$ as
triples $(x,y,\xi)$ where $x\in X_{m\beta}$, $y\in X_{n\gamma}$, and
$\xi\in \Ext(\mathcal{R}_{x},\mathcal{S}_{y})$. Then if $(x,y,\xi)$
and $(x',y',\xi')$ are points of $E'$ that determine isomorphic
representations, it follows that $(x,y)$ and $(x',y')$ lie in the same
orbit for $PGl_{m}\times PGl_{n}$ and therefore may be taken to be
equal and $\xi$ and $\xi'$ determine isomorphic extensions of
$\mathcal{R}_{x}$ on $\mathcal{S}_{y}$. However, the stabiliser in
$Gl_{m}$ of $x$ is the group of automorphisms of $\mathcal{R}_{x}$ and
similarly the stabiliser in $Gl_{n}$ of $y$ is the group of
automorphisms of $\mathcal{S}_{y}$ and so $\xi$ and $\xi'$ are in the
same orbit under the action of the stabiliser of $(x,y)$ and hence
$(x,y,\xi)$ and $(x',y',\xi')$ must lie in the same orbit for
$PGl_{m,n}$ as required.

Finally, it is clear that for each vertex $v$, $\mathcal{T}_{v}$ is a
vector bundle of weight $1$ since it has a subbundle isomorphic to
$\mathcal{S}_{v}$ such that the factor bundle is isomorphic to
$\mathcal{R}_{v}$. 

Thus we have constructed a $(Gl_{m}\times Gl_{n},PGl_{m,n})$-standard
family of representations of dimension vector $\alpha$ for the quiver
$Q$ which is $PGl_{m,n}$ birational to $R(Q(p,t,q),(m\ n))$.

Now we assume that the dimension vector $(m\ n)$ for the quiver
$Q(p,t,q)$ is reducible to matrix normal form. Thus we have a family
$\mathcal{R}$ of representations of dimension vector $(m\ n)$ for the
quiver $Q(p,t,q)$ over the algebraic variety $X$ which is a
$(Gl_{h},PGl_{h})$ standard family such that $X$ is $PGl_{h}$
birational to $M_{h}(k)^{P}$ where $h=\hcf(m,n)=g(\alpha)$ and
$P=1-\langle\alpha/h,\alpha/h\rangle 
=1-\langle(m'\ n'),(m'\ n')\rangle$
where $(m'\ n')=(m\ n)/h$. By lemma \ref{lem:tg}, we may assume that 
$\mathcal{R}(v)\cong k^{m'}\otimes k^{h}\times X$ where $Gl_{h}$ acts
trivially on $k^{m'}$ and diagonally on the remaining terms and
similarly, $\mathcal{R}(w)\cong k^{n'}\otimes k^{h}\times X$ with a
similar action of $Gl_{h}$. By choosing a basis of $k^{m}\cong
k^{m'}\otimes k^{h}$ and of $k^{n}\cong k^{n'}\otimes k^{h}$ we obtain
a morphism from $X$ to $R(Q(p,t,q),(m\ n))$ such that $\mathcal{R}$ is
the pullback of the standard family on $R(Q(p,t,q),(m\ n))$. Moreover,
if we regard $Gl_{h}$ as acting on $k^{n}\cong k^{n'}\otimes k^{h}$
and on $k^{m}\cong k^{m'}\otimes k^{h}$ via its action on $k^{h}$, and
hence on $R(Q(p,t,q),(m\ n))$ and the standard family, we see that
this morphism from $X$ to $R(Q(p,t,q),(m\ n))$ is $PGl_{h}$
equivariant and the morphism between families is $Gl_{h}$
equivariant. 

If we regard $E'$ as a subvariety of $R(Q(p,t,q),(m\ n))$ with the
family $\mathcal{T}$ of representations of dimension vector $\alpha$
for the quiver $Q$ then the pullback of $\mathcal{T}$ to an open
subvariety of $X$ along the $PGl_{h}$ equivariant morphism to
$R(Q(p,t,q),(m\ n))$ is quickly checked to be a $(Gl_{h},PGl_{h})$
standard family of representations of dimension vector $\alpha$ and by
construction $X$ is $PGl_{h}$ birational to $M_{h}(k)^{P}$ which
completes the proof of this theorem.
\end{proof}

Because of this theorem, it is clear inductively that in order to be
able to prove the main theorem, we need to prove it only for quivers
of type $Q(p,t,q)$; this is the content of our next lemma which thus
concludes most of our proof.

\begin{theorem}
\label{th:2vert}
Let $(m\ n)$ be a Schur root for the quiver $Q(p,t,q)$. Then it is
reducible to matrix normal form. 
\end{theorem}
\begin{proof}
If the dimension vector $(m\ n)$ is a Schur root for the quiver
$Q(p',t',q')$ where $(p',t',q')<(p,t,q)$ then the result follows by
induction from lemma \ref{addanarrow}.  This means that we have
certain minimal cases to look at. The results of section \ref{GKQ}
deal with the case where both $p$ and $q$ are zero so we may assume
that at least one of them is non-zero. Thus the remaining minimal
cases are when $(p,t,q)$ is one of $(1,1,1)$, $(p,t,0)$, and $(0,t,q)$
since if both $p$ and $q$ are positive then it is clear that $(m\ n)$
is a Schur root. The first of these in fact reduces to one of the
other two cases since if $n\geq m$ then $(m\ n)$ is a Schur root for
the quiver $Q(0,1,1)$ because (for example) the dimension vector
$(m,n,n)$ for the quiver with vertices $u$, $v$ and $w$ with $1$ arrow
from $u$ to $v$ and $2$ arrows from $v$ to $w$ lies in the fundamental
region for the action of the Weyl group and if $n\leq m$ then we may
apply lemma \ref{lem:dual}. The second and third case are essentially
equivalent by reversing all the arrows of the quiver and applying
lemma \ref{lem:dual} so we shall deal with the case $(0,t,q)$. If $(m\
n)$ is a Schur root for $Q(0,t,q)$ then $m\leq nt$. We shall show that
the dimension vector $(m\ n)$ for the quiver $Q(0,t,q)$ such that
$m\leq nt$ is always a Schur root and determine its moduli space.

If $m=nt$ then the reflection functor at the first vertex together
with lemma \ref{lem:reflection} shows that the moduli space of
representations of dimension vector $(m\ n)$ is reducible to the
moduli space of representations of dimension vector $(0\ n)$ for
$Q(0,t,q)$ which is just $q$ $n$ by $n$ matrices up to simultaneous
conjugacy. If $m<nt$ we may still apply the reflection functor at the
first vertex which preserves the value of $\hcf(m,n)$ followed by
duality to ensure that $0<m\leq\tfrac{nt}{2}$. If $t>2$ or if $t=2$
and $m<n$, then $m<n(t-1)$ and we may conclude by induction that $(m\
n)$ is a Schur root for $Q(0,t-1,q)$ and that it is reducible to
matrix normal form for the quiver $Q(0,t-1,q)$ and hence for the
quiver $Q(0,t,q)$. Thus we may reduce either to the case $(n,n)$ for
$Q(0,2,q)$ which is reducible to $q+1$ $n$ by $n$ matrices up
to simultaneous conjugacy and thus a Schur root with the correct
moduli space; or else we reduce to $(m\ n)$ where $\tfrac{n}{2}\geq
m>0$ for the quiver $Q(0,1,q)$ and it is enough to deal with the case
where $q=1$ since this is a Schur root as we are about to see.

Consider the quiver $Q'$ with $3$ vertices $u$, $v$ and $w$ with $1$
arrow from $u$ to $v$ and $2$ arrows from $v$ to $w$ and the dimension
vector $(m,n,n)$ where $m\leq\tfrac{n}{2}$; this is also clearly a
Schur root since it lies in the fundamental region for the action of
the Weyl group. Also a general representation of this dimension vector
inverts the first arrow from $v$ to $w$ and hence the moduli space of
representations of this dimension vector is birational to the moduli
space of representations of dimension vector $(m\ n)$ for
$Q(0,1,1)$. Now $(m\ m\ 2m)$ is a uniform rigid sub-dimension vector of
$(m\ n\ n)$ and its root is $(1\ 1\ 2)$ which is a real Schur root. Let
$(a\ b\ c)$ be the rigid sub-dimension vector of $(m\ n\ n)$ that is
constructed from this one such that both it and $(m-a\ n-b\ n-c)$ are
uniform in theorem \ref{t6}. By induction on the pair $(m\ n)$ it
follows that the moduli spaces of representations of dimension vector
the root of $(a\ b\ c)$ and the root of $(m-a\ n-b\ n-c)$ satisfy the
theorem since their construction from smaller Schur roots must involve
a smaller pair than $(m\ n)$ and if $m'=\hcf(a,b,c)$ and
$n'=\hcf(m-a,n-b,n-c)$ then $(m'\ n')<(m\ n)$ so again the proof is
complete by induction.
\end{proof}

It remains to show that a dimension vector for a quiver which may have
loops is also reducible to matrix normal form. Let $Q$ be such a
quiver; we recall the double of $Q$ from the last paragraph of section
\ref{internal} whose notation we shall continue to use in the course
of the next proof.

\begin{theorem}
\label{th:complete}
Every Schur root for a quiver is reducible to matrix normal form.
\end{theorem}
\begin{proof}
Let $\alpha$ be a Schur root for the quiver $Q$ and let $\alpha'$ be
the corresponding dimension vector for the double of $Q$, $Q'$. We
have seen that a Schur root for a quiver without loops is reducible to
matrix normal form and therefore we know that $\alpha'$ is reducible
to matrix normal form. If $h=g(\alpha)$, then $h=g(\alpha')$, so we
assume that we have a $(Gl_{h},PGl_{h})$ standard family of
representations of dimension vector $\alpha'$ of the quiver $Q'$,
$\mathcal{R}$, over the affine algebraic variety $X$ such that $X$ is
$PGl_{h})$ equivariantly birational to $M_{h}(k)^{p}$ for some integer
$p$. Let $X'$ be the affine open $PGl_{h}$ equivariant subvariety
where $\mathcal{R}(v,2)$ is invertible. Then $X'$ carries a family of
representations of the quiver $Q$ of dimension vector $\alpha$ by
defining $\mathcal{S}(v)=\mathcal{R}(v,0)$ and
$\mathcal{S}(a)=\mathcal{R}(a,0)\mathcal{R}(v,2)^{-1}$. Using the fact
that the category of representations of the quiver $Q'$ such that
$R(v,2)$ is invertible is equivalent to the category of
representations of the quiver $Q$, it is a simple matter to check that
$\mathcal{S}$ is a $(Gl_{h},PGl_{h})$ standard family and by
construction it is $PGl_{h}$ equivariantly birational to
$M_{h}(k)^{p}$ as required.
\end{proof}
 
It is perhaps worth stating the rationality results that follow from
our main theorem.

\begin{theorem}
\label{rat}
Let $\alpha$ be a Schur root for the quiver $Q$. Then if
$\hcf_{v}(\alpha(v))=n$ where $n=1,2,3 \text{ or }4$, a moduli space
of representations of dimension vector is rational. If $4<n$ and $n$
divides $420$ then a moduli space of representations of dimension
vector $\alpha$ is stably rational. If $n$ is square-free then a
moduli space of representations of dimension vector is retract
rational. 
\end{theorem}
\begin{proof}
This follows from the known results on matrices up to simultaneous
conjugacy. A good summary of the known results may be found in
\cite{leB}.  
\end{proof}

\end{document}